# NIELSEN METHODS AND GROUPS ACTING ON HYPERBOLIC SPACES

ILYA KAPOVICH AND RICHARD WEIDMANN

ABSTRACT. We show that for any $n \in \mathbb{N}$ there exists a constant $C(n)$ such that any $n$-generated group $G$ which acts by isometries on a $\delta$-hyperbolic space (with $\delta > 0$) is either free or has a nontrivial element with translation length at most $\delta C(n)$.

## 1. INTRODUCTION

In his seminal paper on hyperbolic groups M. Gromov states the following:

**Theorem [5.3.A [10]]** *Let $G$ be a $\delta$-hyperbolic group. Suppose that $H$ is a subgroup generated by $n$ elements which contains no non-trivial element of conjugacy length less than $1000\delta \log(n+100)$. Then $H$ is free and quasiconvex in $G$.*

M. Gromov sketches a justification of this statement when $G$ is the fundamental group of a compact negatively curved manifold which relies on edge-angle inequality considerations in CAT(-1) geometry. However, no written proof of the above theorem exists in the literature so far.

In this paper we will prove the following natural generalization of Gromov's statement:

**Theorem 1.** *For any integer $n > 0$ there exists a constant $C(n)$ with the following property.*

*Suppose a group $G$ generated by elements $g_1, \ldots, g_n$ acts by isometries on a $\delta$-hyperbolic geodesic metric space $(X, d)$, where $\delta > 0$. Then one of the following holds:*
*(1) The group $G$ is free on $M = (g_1, \ldots, g_n)$ and for any $x \in X$ the map $G \to X$ defined by $h \mapsto hx$ is a quasi-isometric embedding (where $G$ is considered with the free group word-metric).*
*(2) The tuple $M$ is Nielsen-equivalent to a tuple $M' = (f_1, \ldots, f_n)$ such that $d(f_1 y, y) < \delta C(n)$ for some $y \in X$.*

That is to say, a group acting by isometries on a hyperbolic space is either free or contains a nontrivial element of small translation length. Note that the statement of Theorem 1 is not entirely obvious even if the action of $G$ on $X$ is non-discrete or non-faithful or has elliptic or parabolic elements. In each of these cases it is clear that $G$ has an element of small translation length, but it is not obvious that such an element can be made a member
of the generating set of $G$ with minimal cardinality. It is easy to see that if
case (1) of Theorem 1 holds then the orbit $Gx$ is a quasiconvex and discrete
subset of $X$. Moreover, the action of $G$ on $X$ is discontinuous in the sense
that for any bounded subset $K \subset X$ the set $\{g \in G | gK \cap K \neq \emptyset\}$ is finite. In
particular, if the space $X$ is proper (that is all closed balls in $X$ are compact)
then the action of $G$ on $X$ is properly discontinuous. It appears that the
conclusion of Theorem 1 is new even for groups acting on the standard
hyperbolic space $\mathbb{H}^m$.

Recall that $\mathbb{R}$-trees are precisely 0-hyperbolic geodesic metric spaces.
Thus an $\mathbb{R}$-tree is $\delta$-hyperbolic for any $\delta > 0$. By considering arbitrarily
small $\delta > 0$ Theorem 1 implies that an $n$-generated group acting by isome-
tries on an $\mathbb{R}$-tree is either a free group of rank $n$ or after some Nielsen
transformations the first generator can be made to have arbitrarily small
translation length:

**Corollary 2.** *Let* $G = \langle g_1, \ldots, g_n \rangle$ *be an $n$-generated group acting by isome-
tries on an $\mathbb{R}$-tree $X$.*

*Then either $G$ is a free group of rank $n$ with free basis $g_1, \ldots, g_n$ which
acts freely on $X$ with the orbit map being a quasi-isometric embedding or
for any $\epsilon > 0$ there is a tuple $(h_1, \ldots, h_n)$ Nielsen-equivalent to $(g_1, \ldots, g_n)$
such that $||h_1|| \leq \epsilon$.*

*(Here for $g \in G$ $||g|| := \inf\{d(x, gx) | x \in X\}$ is the translation length of
$g$).*

Again, the above fact is not obvious even if $G$ acts non-freely or non-
discretely on an $\mathbb{R}$-tree $X$, although it can probably be obtained using Rips'
theory [3, 12, 15].

A word-hyperbolic group $G$ has a finite generating set $X$ such that the
Cayley graph $\Gamma(G, X)$ with respect to $X$ is $\delta$-hyperbolic. The group $G$ has
a canonical isometric action on $\Gamma(G, X)$. For this action if $y \in G$ is a vertex
of $\Gamma(G, X)$, $g \in G$ and $d(y, gy) = c$ then the element $y^{-1}gy \in G$ has length
$c$ in the word metric with respect to $X$. Hence Theorem 1 immediately
implies Gromov's claim about subgroups of hyperbolic groups, except that
our result does not yield the same constant. Unlike M. Gromov, we do not
claim anything about the asymptotics of $C(n)$ in terms of $n$. Nevertheless, a
careful analysis of our argument shows that our constant $C(n)$ in Theorem 1
is a recursive function of $n$ and thus can be computed algorithmically.

T. Delzant proved in [6] that a torsion-free word-hyperbolic group has only
finitely many conjugacy classes of non-free two-generated subgroups. In this
paper he also proves Theorem 1 for $n = 2$. M. Koubi [19] has established
Theorem 1 for $n = 3$. We are informed that Goulnara Arzhantseva [1]
independently of our work proved Theorem 1 for hyperbolic groups using a
method different from the one employed in the present paper.

Our strategy for obtaining the main result relies on the use of Nielsen
methods. The theory of Nielsen transformations was introduced by Jacob

Nielsen in the 1920s (c.f [21]) to show that subgroups of free groups are free. H. Zieschang [29] developed an analog of Nielsen's theory to study subgroups of amalgamated free products. A similar theory has been developed for more general splittings by N. Peczynski and W. Reiwer [25], S. Pride [22], A.Hoare [14], and the authors of the present paper in [16] and [28]. Nielsen methods proved to be particularly useful for studying discreteness properties for groups generated by several matrices in $PSL(2, \mathbb{R})$ and $PSL(2, \mathbb{C})$ (see for example [11, 8, 9]). Various other applications of Nielsen methods can be also found in [23], [5], [14], [27], [24], [7], [17] and other sources.

The idea of Nielsen's approach is to study the subgroup generated by an $n$-tuple of elements by applying a chain of simple moves, called Nielsen transformations, which change the $n$-tuple but preserve the subgroup it generates. Then one looks at a minimal (in some natural sense) $n$-tuple which is Nielsen-equivalent to the original $n$-tuple and one tries to analyze the minimal case to obtain some information about the subgroup generated by this tuple. It turns out that all truly powerful generalizations of Nielsen method (e.g the work of H. Zieschang) employ the notion of minimality which incorporates a combinatorial as well as a geometric component. In fact it is the combinatorial component that usually gives the method its full power. Thus both J. Nielsen and H. Zieschang use some sort of lexicographical ordering on the group. The main technical difficulty we had to overcome was to find a suitable substitution for the lexicographical ordering (or, rather, the effect that this ordering should produce) which would work for groups acting on hyperbolic spaces.

It is worth noting that in the proof of the theorem regarding two-generator subgroups of hyperbolic groups mentioned above [6] T. Delzant also essentially uses Nielsen methods. However, the argument needed for the case $n = 2$ is considerably less complicate as the reduction process uses simple Nielsen moves while in the case $n \geq 3$ we sometimes apply multiple Nielsen moves at once. This is due to the fact that the cancellation in much longer products has to be investigated.

The statement of Theorem 1, while not giving an explicit asymptotics of $C(n)$, is in other regards considerably stronger than the original claim of M. Gromov. Indeed, the fact that we can bound the length of the first generator (and not just of some element) of the group $H$ can, under appropriate assumptions, be considerably generalized and pushed much further. Thus, in [18] we prove close analogs of the results of [28]. One application is the following: If $G$ is torsion-free word-hyperbolic group where all $k$-generated subgroups are quasiconvex, then $G$ contains only finitely many conjugacy classes of one-ended $(k + 1)$-generated subgroups. Moreover, in such $G$ for any $n \geq k + 1$ any $n$-tuple of elements generating an one-ended subgroup is Nielsen-equivalent (after conjugation) to an $n$-tuple where the first $k+1$ elements are short. Since 1-generated (i.e. cyclic) subgroups are always quasiconvex in a hyperbolic group, this implies the result of T. Delzant mentioned above. Other applications of the methods developed in this paper and proved in [18] include the acylindrical accessibility theorem for finitely generated groups acting on $\mathbb{R}$-trees.

## 2. Preliminaries

We recall the definition of Nielsen equivalence and some basic facts about hyperbolic spaces although we assume some familiarity of the reader with this theory. Background information can be found in [10], [2], [4], [13].

**Definition 3** (Nielsen equivalence). Let $G$ be a group and let $M = (g_1, \ldots, g_n) \in G^n$ be an $n$-tuple of elements of $G$. The following moves are called *elementary Nielsen moves* on $M$:

(N1) For some $i, 1 \leq i \leq n$ replace $g_i$ by $g_i^{-1}$ in $M$.
(N2) For some $i \neq j$, $1 \leq i, j \leq n$ replace $g_i$ by $g_i g_j$ in $M$.
(N3) For some $i \neq j$, $1 \leq i, j \leq n$ interchange $g_i$ and $g_j$ in $M$.

We say that $M = (g_1, \ldots, g_n) \in G^n$ and $M' = (f_1, \ldots, f_n) \in G^n$ are *Nielsen-equivalent*, denoted $M \sim_N M'$, if there is a chain of elementary Nielsen moves which transforms $M$ to $M'$.

It is easy to see that Nielsen equivalence is indeed an equivalence relation on $G^n$. Moreover, if $M \sim_N M'$ then $\langle M \rangle = \langle M' \rangle \leq G$, that is Nielsen-equivalent tuples generate the same subgroup in $G$. The above definition also implies that the following move preserves the Nielsen-equivalence class of $M = (g_1, \ldots, g_n) \in G^n$:

(N4) For some $i, 1 \leq i \leq n$ replace $g_i$ by $h g_i^{\epsilon} h'$ where $\epsilon \in \{-1, 1\}$, $h, h' \in \langle g_1, \ldots, g_{i-1}, g_{i+1}, \ldots, g_n \rangle \leq G$.

For any two points $x, y$ in a geodesic metric space $(X, d)$ we will often denote by $[x, y]$ a geodesic segment from $x$ to $y$.

**Definition 4** (Hyperbolic space). We say that a geodesic metric space $(X, d)$ is $\delta$-*hyperbolic*, provided that all geodesic triangles are $\delta$-thin, that is if the following holds:

For any $x, y, z \in X$ and geodesic segments $[x, y]$, $[x, z]$ and $[y, z]$ the segment $[x, y]$ is contained in the $\delta$-neighborhood of $[x, z] \cup [y, z]$.

Since we will only need to deal with 1-hyperbolic spaces, we will formulate the following facts for 1-hyperbolic spaces only.

**Definition 5** (Quasigeodesics and local quasigeodesics). Let $(X, d)$ be a metric space. Let $\alpha : I \longrightarrow X$ be a naturally parameterized path in $X$ (here $I$ is a subinterval of the real line). Let $\lambda > 0, \epsilon \geq 0, T > 0$ be some numbers.

(1) The path $\alpha$ is called a $(\lambda, \epsilon)$-*quasigeodesic* in $X$ if for any closed subinterval $[a, b] \subseteq I$ we have

$$|b - a| \leq \lambda d(\alpha(a), \alpha(b)) + \epsilon.$$

(2) The path $\alpha$ is called a $T$-*local* $(\lambda, \epsilon)$-*quasigeodesic* in $X$ if for any closed subinterval $[a, b] \subseteq I$ with $|b - a| \leq T$ the restriction of $\alpha$ to $[a, b]$ is a $(\lambda, \epsilon)$-quasigeodesic

It is well-known that in hyperbolic metric spaces quasigeodesics are Hausdorff close to geodesics and that local quasigeodesics are global quasigeodesics, provided the local parameter $T$ is chosen big enough. Namely, the following holds:

**Lemma 6.** *For any $\lambda > 0$, $\epsilon \geq 0$ there exist constants $H = H(\lambda, \epsilon) > 0$, $T = T(\lambda, \epsilon) > 0$ and $K = K(\lambda, \epsilon) > 0$ such that the following hold:*
*(a) Let $\alpha : I \longrightarrow X$ be a $(\lambda, \epsilon)$-quasigeodesic in a 1-hyperbolic metric space $(X, d)$. Then for any $[a, b] \subset I$ and any geodesic $\gamma = [\alpha(a), \alpha(b)]$ in $X$ the path $\alpha|_{[a,b]}$ and $\gamma$ are $H$-Hausdorff close.*
*(b) Any $T$-local $(\lambda, epsilon)$-quasigeodesic in a 1-hyperbolic metric space $(X, d)$ is $(K, K)$-quasigeodesic.*

*Proof.* Part (a) is very well-known and is proved, for example in Ch. 5, Theorem 11 of [13]. Part (b) is equally well-known (and sometimes referred to as "pasting quasigeodesics") and is proved in [4]. It also follows immediately from Ch. 5, Theorem 21 of [13]                                       □

**Convention 7** (Constants)**.** Unless otherwise stated, we will assume that all constants in this paper are positive integers. There exist integers $N_0 > 0, K > 0$ with the following property: Any $N_0$-local $(1, 100)$-quasigeodesic in a 1-hyperbolic space is $(K, K)$-quasigeodesic. Let $L > 0$ be an integer such that any two $(K, K)$-quasigeodesics in a 1-hyperbolic space are $L$-Hausdorff close. The existence of such $K, L, N_0$ follows from Lemma 6. Define $N_1 := 1000 L K^2 N_0$. The constants $K, L, N_0, N_1$ will be fixed for the remainder of this paper.

If $S$ is a subset of a group $G$ we will denote by $\langle S \rangle$ the subgroup of $G$ generated by $S$.

## 3. Main technical results

Till the remainder of this paper, unless stated otherwise, let $X$ be a 1-hyperbolic geodesic metric space with a base-point $x \in X$ and let $G$ be a group acting on $X$ by isometries.

For any element $g \in G$ we denote by $\|g\|$ the *translation length* of $g$ defined as

$$\|g\| := \inf_{y \in X} d(y, gy).$$

We further define $|g|_y = d(y, gy)$ for $y \in X$, $g \in G$.

For a $n$-tuple $M = (g_1, \dots, g_n) \in G^n$ we put $|M|_x := \sum_{i=1}^n |g_i|_x$. If the choice of base-point $x$ in unambiguous, we will often omit the subscript and use the notation $|M| = |M|_x$.

**Definition 8.** Minimal tuple say that an $n$-tuple $M \in G^n$ is *minimal* if

$$|M| \leq |M'| + 1$$

for any $n$-tuple $M' \in G^n$ which is Nielsen-equivalent to $M$. Thus

$$|M| \leq 1 + \inf\{|M'| \text{ where } M' \sim_N M\}.$$

**Remark 9.** It is easy to see that any set $M$ is Nielsen-equivalent to a minimal set since the infimum in the above formula can be approximated with arbitrarily small error.

**Remark 10** (An informal remark about the Nielsen method). It is worth-while to point out the crucial role played by the lexicographical order in the original approach of Nielsen for the free group case. We will sketch a rough outline of Nielesen's argument in this case.

Namely, suppose $F = F(X)$ is a free group on $X = x_1, \ldots x_m$ and that we want to study the subgroup $H$ of $F$ generated by elements $h_1, \ldots, h_n \in F$. Fix a linear order of the set $X \cup X^{-1}$ and the induced lexicographical order on the set of words in the alphabet $X \cup X^{-1}$.

J. Nielsen proceeded by performing elementary Nielsen moves on the $n$-tuple $M = (h_1, \ldots, h_n)$ to minimize the sum $|M|$ of the lengths of its components $h_i$ (considered as freely reduced words in $X$). This is done step-by-step by investigating cancellation in the products of the type $h_i h_j^{\pm 1}$. For example, suppose that more than a half of $h_j$ cancels in the product $h_i h_j^{\epsilon}$, where $\epsilon \in \{-1, 1\}$ and $i \neq j$. Then $|h_i h_j^{\pm 1}| < |h_i|$ and we replace $h_i$ by $h_i h_j^{\epsilon}$ while leaving other members of the tuple $(h_1, \ldots, h_n)$ intact. This is clearly an elementary Nielsen move producing a new $n$-tuple $M'$ with $|M'| < |M|$.

It may happen, however, that exactly a half of $h_j$ cancels in the product $h_i h_j^{\epsilon}$. In this case $|h_i h_j^{\epsilon}| = |h_i|$ and replacing $h_i$ by $h_i h_j^{\epsilon}$ preserves $|M|$. Thus it is not clear whether and when such a move should be performed. Lexicographical ordering solves this difficulty. We have $h_i = uv^{-1}$, $h_j^{\epsilon} = vw^{-1}$, where $|v| = |w| = |h_j|/2$ and $v$ is the front half of $h_j^{\epsilon}$ that cancels in the product $h_i h_j^{\epsilon}$. Then $h_i h_j^{\epsilon} = uw^{-1}$. Roughly speaking, we perform the Nielsen move $h_i \mapsto h_i h_j^{\epsilon}$ only if $w$ is lexicographically smaller than $v$.

This process (which has to be defined more precisely) is performed repeatedly. Lexicographical ordering guarantees that even if a particular step preserves the length $|M|$, some finer type of complexity decreases at each step. It may happen that after a number of steps we obtain a tuple where some entry is equal to 1. In this case the entry is crossed-out and the whole process is applied again to the resulting $(n-1)$-tuple. Eventually the process does terminate with a tuple of elements where all entries are non-trivial and no step of the above types can be applied. This tuple is then shown to be a free basis of $H$.

Several difficulties arise when this kind of argument is "quasified" and applied to a group acting on a hyperbolic space. First, lexicographical order is not available and some other trick is needed in its place. Second, much more care is needed in order to define what it means for "approximately a half" of $h_j$ to cancel in the product $h_i h_j^{\epsilon}$. Specifically, any such definition should involve "approximation" constants with inductive dependence on $n$. We solve these problems by proving the following theorem and using it as a main technical tool of the paper:

**Theorem 11.** *Let $n \in N$ and $c$ be positive numbers. Then there exist numbers $d_1 = d_1(n, c)$, $d_2 = d_2(n)$, $d_3 = d_3(n)$, and $d_4 = d_4(n, c)$ such that every minimal $n$-tuple $M = (g_1, \ldots, g_n) \in G^n$ is either Nielsen-equivalent to an $n$-tuple $M' = (g_1', \ldots, g_n')$ such that $||g_1'|| \leq d_1$ or $U = \langle M \rangle \leq G$ is freely generated by $M$ and the following hold:*

1. *The map $U \to X$ defined by $u \mapsto ux$ is a quasi-isometric embedding, where $U$ is considered with the free group word-metric.*
2. *For any $u = g_{i_1}^{\epsilon_1} \cdots g_{i_k}^{\epsilon_k} \in U$, the segment $[x, ux]$ lies in the $d_2$-neighborhood of $[x, g_{i_1}^{\epsilon_1} x] \cup g_{i_1}^{\epsilon_1}[x, g_{i_2}^{\epsilon_2} x] \cup \cdots \cup g_{i_1}^{\epsilon_1} \cdots g_{i_{k-1}}^{\epsilon_{k-1}}[x, g_{i_k}^{\epsilon_k} x]$. This implies in particular that $[x, ux]$ is contained in the $a$-neighborhood of $Ux$ where $a = \max\limits_{i=1,\ldots,n}(|g_i|_x/2 + d_2)$.*
3. *For any freely reduced $u = g_{i_1}^{\epsilon_1} \cdots g_{i_k}^{\epsilon_k} \in U$ we have that $|u|_x \geq |g_{i_j}|_x - d_3$.*
4. *If a $S$ is a subsegment of geodesic segment $[x, ux]$ for some $u \in U$ where the length of $S$ is greater than $d_4$ then $S$ intersects nontrivially the $b$-neighborhood of $Ux$ where $b = \max\limits_{i=1,\ldots,n}(|g_i|_x/2 - c)$.*

Note that if $x' \in X$ then $d(ux, ux') = d(x, x')$ and $|d(x, ux) - d(x', ux')| \leq 2d(x, x')$ for any $u \in U$. Hence the orbit maps of $U$ with respect to $x$ and $x'$ are $2d(x, x')$-close to each other. Thus if the $U$-orbit map with respect to $x$ is a quasi-isometric embedding, then so is the orbit map of $U$ with respect to any other point $x' \in X$.

The above theorem now immediately implies Theorem 1.

*Proof of Theorem 1*

We put $C(n) := d_1(n, 1)$, where $d_1$ is the constant provided by Theorem 11.

If $(X, d)$ is a 1-hyperbolic space, then the conclusion of Theorem 1 follows directly from Theorem 11.

Suppose now $(X, d)$ is a $\delta$-hyperbolic space for some $\delta > 0$. Then the scaled version of $X$, namely $(X, d/\delta)$ is 1-hyperbolic. Applying the result already known for 1-hyperbolic spaces we immediately obtain the general statement of Theorem 1. □

The advantage of Theorem 11 over Theorem 1 is that we are able to prove it by induction on $n$. This indeed will be our strategy.

To establish the base of induction we need the following simple lemma.

**Lemma 12.** *There exists a constant $c$ such that for any element $g \in G$ one of the following holds:*

1. *$||g|| \leq c$;*
2. *there exists a point $y \in X$ such that for any $n \in \mathbb{N}$ the path*

$$[x, y] \cup [y, gy] \cup g[y, gy] \cup \cdots \cup g^{n-1}[y, gy] \cup [g^n y, g^n x]$$

*is a $(c, c)$-quasigeodesic lying in the $c$-neighborhood of any geodesic segment $[x, g^n x]$.*

*Proof.* We will sketch the proof of this lemma and leave the details to the reader. Let $z \in X$ be such that $d(z, gz) \leq ||g|| + 1$, so that $y$ "almost" realizes the translation length of $g$. Denote $N = d(z, gz)$. If $N = 0$ then part (1) of the lemma obviously holds. Suppose $N > 0$.

Consider the biinfinite path

$$\sigma = \cdots \cup g^{-1}[z, gz] \cup [z, gz] \cup g[z, gz] \cup g^2[z, gz] \cup \ldots$$

By the choice of $z$ the path $\sigma$ is $N$-local $(1,2)$-quasigeodesic. Therefore there is a constant $c'$ such that if $N \geq c'$ then $\sigma$ is a global $c'$-quasigeodesic and for any two points $a, b$ on $\sigma$ the segment of $\sigma$ between $a$ and $b$ is $c'$-Hausdorff close to any geodesic segment $[a, b]$. Put $y$ to be the "nearest point projection" of $x$ onto $\sigma$. That is let $y \in \sigma$ be such that $d(x, y) \leq d(x, y') + 1$ for any $y' \in \sigma$. It is now easy to see that this $y$, with an appropriate $c$, satisfies condition (2) in the lemma. $\square$

The key to the proof of Theorem 11 is the following proposition.

**Proposition 13.** *Let $T \geq N_1$ be a constant and suppose that Theorem 11 holds for $n - 1$. Then there exist numbers $c_1 = c_1(n, T)$, $c_2 = c_2(n)$ and $c_3 = c_3(n)$ such that any minimal $n$-tuple $M = (g_1, \ldots, g_n) \in G^n$ with $|g_i|_x \leq |g_n|_x$ for $1 \leq i \leq n - 1$ is either Nielsen-equivalent to a tuple $M' = (g'_1, \ldots, g'_n)$ such that $||g'_1|| \leq c_1$ or $U = \langle M \rangle$ is freely generated by $M$ and the following hold:*

1. *The map $U \to X$ defined by $u \mapsto ux$ is a quasi-isometric embedding (where $U$ is considered with the free group word-metric).*
2. *For any $u = g_{i_1}^{\epsilon_1} \cdots g_{i_k}^{\epsilon_k} \in U$, the segment $[x, ux]$ lies in the $c_2$-neighborhood of $[x, g_{i_1}^{\epsilon_1} x] \cup g_{i_1}^{\epsilon_1}[x, g_{i_2}^{\epsilon_2} x] \cup \cdots \cup g_{i_1}^{\epsilon_1} \cdots g_{i_{k-1}}^{\epsilon_{k-1}}[x, g_{i_k}^{\epsilon_k} x]$.*
3. *For any freely reduced $u = g_{i_1}^{\epsilon_1} \cdots g_{i_k}^{\epsilon_k} \in U$ we have that $|u|_x \geq |g_{i_j}|_x - c_3$.*
4. *For any product of type $w = h_1 g_n^{\epsilon_1} h_2 g_n^{\epsilon_2} \cdots g_n^{\epsilon_{l-1}} h_l g_n^{\epsilon_l} h_{l+1}$ with $\epsilon_i \in \{-1, 1\}$ for $1 \leq i \leq l$, $h_i \in \langle g_1, \ldots, g_{n-1} \rangle$ for $1 \leq i \leq l + 1$ and $h_i \neq 1$ if $\epsilon_i = -\epsilon_{i+1}$ for $2 \leq i \leq l$ the following holds:*
   *The $2L + 100$-neighborhood of every subsegment of $[x, wx]$ of length at least $10T$ contains a subsegment of length at least $T$ of either a segment of type $h_1 g_n^{\epsilon_1} h_2 g_n^{\epsilon_2} \cdots g_n^{\epsilon_{i-1}} h_i[x, g_n^{\epsilon_i} x]$ or $h_1 g_n^{\epsilon_1} h_2 g_n^{\epsilon_2} \cdots g_n^{\epsilon_{i-1}}[x, h_i x]$ for some $i$.*

We first show how to deduce Theorem 11 from the Proposition 13 and then devote the remainder of this paper to the proof of Proposition 13.

*Proof of Theorem 11.* The proof is by induction on $n$. In the case $n = 1$ the existence of the constants $d_1, \ldots, d_4$ easily follows from Lemma 12.

Suppose now that $n > 1$ and that Theorem 11 has been established for $n - 1$. Put

$$T := \max(N_1, d_4(n - 1, c + 2L + 100 + 1), 2c + 2(2L + 100) + 1)$$

and put $d_1(n, c) := c_1(n, T)$, where $c_1$ is the constant provided by Proposition 13. Put $d_2(n) := c_2(n)$ and $d_3(n) := c_3(n)$. Finally, put $d_4(n, c) := \max\{10T, d_4(n - 1, c)\}$.

Let $M = (g_1, \ldots, g_n) \in G^n$ be a minimal $n$-tuple where $|g_1|_x \leq |g_2|_x \leq \cdots \leq |g_n|_x$. Suppose that $|g_1|_x > d_1(n, c)$. It is clear that assertion 1-3 of Theorem 11 follow directly from Proposition 13.

Thus it suffices to prove that assertion (4) of Theorem 11 holds. In the following we denote the subgroup $\langle g_1, \ldots, g_{n-1} \rangle$ of $G$ by $H$.

Suppose now that $u \in U$ and that $S$ is a subsegment of a geodesic $[x, ux]$ such that the length of $S$ is greater than $d_4(n, c)$. If $u \in H$, then part (4) of Theorem 11 holds by the inductive hypothesis since the length of $S$ is at least $d_4(n - 1, c)$.

If $u \in U - H$ then $u$ can be represented as the product of the form:

$$u = h_1 g_n^{\epsilon_1} h_2 g_n^{\epsilon_2} \cdots g_n^{\epsilon_{l-1}} h_l g_n^{\epsilon_l} h_{l+1}$$

with $\epsilon_i \in \{-1, 1\}$ for $1 \leq i \leq l$, $h_i \in H$ for $1 \leq i \leq l + 1$ and $h_i \neq 1$ if $\epsilon_i = -\epsilon_{i+1}$ for $2 \leq i \leq l$.

Note that $S$ has length at least $10T$ by the choice of $d_4(n, c)$. By Proposition 13 this implies that the $(2L+100)$-neighborhood of $S$ contains a subsegment $SS$ of length at least $T$ of one of the segments $h_1 g_n^{\epsilon_1} h_2 g_n^{\epsilon_2} \cdots g_n^{\epsilon_{i-1}} h_i [x, g_n^{\epsilon_i} x]$ or $h_1 g_n^{\epsilon_1} h_2 g_n^{\epsilon_2} \cdots g_n^{\epsilon_{i-1}} [x, h_i x]$.

If $SS$ is a subsegment of a segment of type $h_1 g_n^{\epsilon_1} h_2 g_n^{\epsilon_2} \cdots g_n^{\epsilon_{i-1}} h_i [x, g_n^{\epsilon_i} x]$ then there is point $q$ of $SS$ that is in distance less that $|g_n|_x/2 - (c + 2L + 100)$ of one of the extremals of the segment $h_1 g_n^{\epsilon_1} h_2 g_n^{\epsilon_2} \cdots g_n^{\epsilon_{i-1}} h_i [x, g_n^{\epsilon_i} x]$ since

$$l(SS) \geq T \geq 2c + 2(2L + 100) + 1.$$

Thus $d(q, u'x) \leq |g_n|_x/2 - (c + 2L + 100)$ for some $u' \in U$. Since $SS$ lies in the $(2L + 100)$-neighborhood of $S$ this implies that there is a point $p$ of $S$ with $d(p, q) \leq 2L + 100$. Therefore

$$d(p, u'x) \leq d(p, q) + d(q, ux') \leq (2L + 100) + |g_n|_x/2 - (c + 2L + 100) = |g_n|_x/2 - c$$

and so assertion (4) of Theorem 11 holds.

Suppose now that $SS$ is a subsegment of type $h_1 g_n^{\epsilon_1} h_2 g_n^{\epsilon_2} \cdots g_n^{\epsilon_{i-1}} [x, h_i x]$. Then there is a point $q$ of $SS$ that is in distance at most

$$b = \max_{i=1,\ldots,n-1} (|g_i|_x/2 - (c + 2L + 100 + 1)) \leq |g_n|_x/2 - (c + 2L + 100 + 1)$$

of some point of $Ux$ since statement (4) of Theorem 11 holds for $n - 1$ and since $T \geq d_4(n - 1, c + 2L + 100 + 1)$. That is for some $u' \in U$ we have $d(q, u'x) \leq |g_n|_x/2 - (c + 2L + 100 + 1)$. Since $SS$ is in the $2L + 100$-neighborhood of $S$, there is a point $p$ on $S$ such that $d(p, q) \leq 2L + 100$.

Therefore again we have

$$d(p, u'x) \leq d(p, q) + d(q, u'x) \leq (2L + 100) + |g_n|_x/2 - (c + 2L + 100) = |g_n|_x/2 - c$$

and so assertion (4) of Theorem 11 holds. This completes the proof of Theorem 11. $\square$

## 4. Auxiliary lemmas

Until the end of the paper, unless specified otherwise, let $G$ be a group generated by elements $g_1, \ldots, g_n$ on a 1-hyperbolic geodesic metric space $X$ with a base-point $x \in X$. We will denote $M = (g_1, \ldots, g_n)$ and $M_{n-1} = (g_1, \ldots, g_{n-1})$. We also denote $H = \langle M_{n-1} \rangle = \langle g_1, \ldots, g_{n-1} \rangle \leq G$.

In this section we will look at products in $g_n^{\pm 1}$ and elements of the subgroup $H$. Our aim is to show that under suitable minimality assumptions on $M$ "local cancellation" in $X$ corresponding to such products is limited.

All statements in this section make the following important assumption:

**Convention 14** (Minimal)**.** Throughout this section we assume that the $n$-tuple $M = (g_1, \ldots, g_n) \in G^n$ is minimal with $|g_i|_x \leq |g_n|_x$ for $1 \leq i \leq n-1$ and suppose that Theorem 11 holds for $n - 1$.

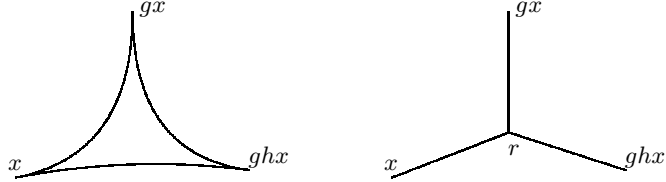

FIGURE 1. Approximating triangles in a 1-hyperbolic space

We first look at products of "syllable length" two and three. The following lemmas follow immediately from the definition of hyperbolicity.

**Lemma 15.** Let $g, h \in G$. Then there is a path $w_g = [x, r] \cup [r, gx]$ of length at most $|g|_x + 2$ and a path $w_h = [gx, r] \cup [r, ghx]$ of length at most $|h|_x + 2$ such that $w_{gh} = [x, r] \cup [r, ghx]$ is a path of length at most $|gh|_x + 2$ (see Figure 1) . Also $w_k$ lies in the 2-neighborhood of any geodesic segment joining their endpoints for $k \in \{g, h, gh\}$.

**Lemma 16.** Let $g, h, f \in G$. Then one of the following occurs:

1. There are paths $w_g = [x, p] \cup [p, gx]$ of length at most $|g|_x + 4$, $w_h = [gx, p] \cup [p, q] \cup [q, ghx]$ of length at most $|h|_x + 4$ and $w_f = [ghx, q] \cup [q, ghfx]$ of length at most $|g|_x + 4$ and the path $w_{ghf} = [x, p] \cup [p, q] \cup [q, ghfx]$ is of length at most $|ghf|_x + 4$ (see Figure 2).

2. There are paths $w_g = [x, p] \cup [p, q] \cup [q, gx]$ of length at most $|g|_x + 4$, $w_h = [gx, q] \cup [q, ghx]$ of length at most $|h|_x + 4$ and $w_f = [ghx, q] \cup [q, p] \cup [q, ghfx]$ of length at most $|g|_x + 4$ and the path $w_{ghf} = [x, p] \cup [p, ghfx]$ is of length at most $|ghf|_x + 4$ and the segment $[p, q]$ lies in the 4-neighborhood of any of the geodesic segments $[x, gx]$ and $[ghx, ghfx]$ (see Figure 3). Furthermore the 5-neighborhood of $[p, q]$ contains a segment of $[x, gx]$ and of $gh[x, fx]$ of length at least $d(p, q)$.

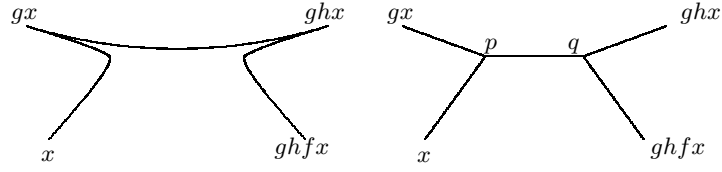

FIGURE 2. Thin quadrilaterals: first type

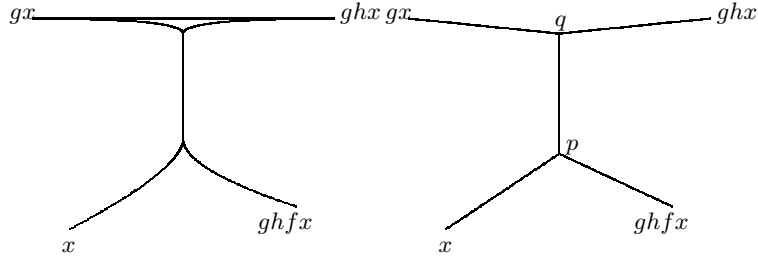

FIGURE 3. Thin quadrilaterals: second type

*Also $w_k$ lies in the 4-neighborhood of any geodesic segment joining their endpoints for $k \in \{g, h, f, ghf\}$ and the paths can be chosen such that for any choice of segments and $r$ as in Lemma 15 we can assume that the following holds: If the first (second) case of Lemma 16 occurs then $[x, r]$ lies in the 5-neighborhood of $[x, p]$ ($[x, p] \cup [p, q]$), $[r, gx]$ lies in the 5-neighborhood of $[p, gx]$ ($[q, gx]$) and $[r, ghx]$ lies in the 5-neighborhood of $[p, q] \cup [q, ghx]$ ($[q, ghx]$).*

The above lemma easily implies the following:

**Corollary 17.** *Let $f, g, h \in G$. Let $\tau$ be a subsegment of length $3D$ in $[x, ghfx]$. Then the 24-neighborhood of $\tau$ contains a subsegment of length at least $D - 40$ of one of $[x, gx]$, $g[x, hx]$, $gh[x, fx]$.*

**Convention 18.** Throughout this section, unless specified otherwise, notation $N$ will stand for an integer $N \geq N_1$.

The following important but simple lemma is an immediate consequence of Lemma 15. It states that, under the minimality assumption on $M$ made in this section, approximately a half of $[x, hx]$ and approximately a half of $h[x, g_n^\epsilon x]$ "survive" in $[x, hg_n^\epsilon x]$ (where $h \in H$, $\epsilon = \pm 1$).

**Lemma 19.** *Let $N \geq N_1$ be an integer. Suppose that $h \in H$, $\epsilon \in \{-1, 1\}$ and that $d_2 = d_2(n-1)$ from the conclusion of Theorem 11. Then the following hold:*

1. *Either $||g_n|| < 2N + 2d_2 + 5$ or there exists a path $w = [x, r] \cup [r, g_n^2 x]$ of length at most $|g_n^2|_x + 2$ in the 2-neighborhood of $[x, g_n^2 x]$ where $[x, r]$ lies*

   *in the 2-neighborhood of $[x, g_n x]$ and $[r, g_n^2 x]$ lies in the 2-neighborhood of $g_n[x, g_n x] = [g_n x, g_n^2 x]$ and $d(x, r) \geq |g_n|_x/2 + N + d_2$ and $d(r, g_n^2 x) \geq |g_n|_x/2 + N + d_2$.*

2. *There exists a path $w = [x, r] \cup [r, hg_n^\epsilon x]$ of length at most $|hg_n^\epsilon|_x + 2$ in the 2-neighborhood of $[x, hg_n^\epsilon x]$ where $[x, r]$ lies in the 2-neighborhood of $[x, hx]$, $[r, hg_n^\epsilon x]$ lies in the 2-neighborhood of $h[x, g_n^\epsilon x] = [hx, hg_n^\epsilon x]$ and $d(r, x) \geq |h|_x/2 - 3$ and $d(r, hg_n^\epsilon x) \geq |g_n|_x/2 - (3 + d_2)$.*

*Proof* The first assertion follows immediately from Lemma 15. The existence of the path $w = [x, r] \cup [r, hg_n^\epsilon x]$ with all required properties but the inequalities follows from Lemma 15. The inequality $d(r, x) \geq |h|_x/2 - 3$ holds since otherwise $|hg_n^\epsilon|_x = d(x, hg_n^\epsilon x) < d(hx, hg_n^\epsilon x) - 1 = d(x, g_n^\epsilon x) - 1 = |g_n|_x - 1$ which contradicts the minimality of $M$.

It remains to show that $d(r, hg_n^\epsilon x) \geq |g_n|_x/2 - (3 + d_2)$. Suppose that $d(r, hg_n^\epsilon x) < |g_n|_x/2 - (3 + d_2)$. Note that $r$ lies in the 2-neighborhood of $[x, hx]$ and that $x, hx \in Hx$. Choose $z \in [x, hx]$ such that $d(r, z) \leq 2$. Since Theorem 11 (2) holds for $n - 1$ we know that $z$ lies in the $a$-neighborhood of $Hx$ where $a = \max_{i=1,\ldots,n-1} (|g_i|_x/2 + d_2) \leq |g_n|_x/2 + d_2$. Choose $\bar{h} \in H$ such that $d(\bar{h}x, z) \leq |g_n|_x/2 + d_2$, it follows that $d(\bar{h}x, r) \leq |g_n|_x/2 + d_2 + 2$. Now this implies that

$$|\bar{h}^{-1}hg_n^\epsilon|_x = d(x, \bar{h}^{-1}hg_n^\epsilon x) = d(\bar{h}x, hg_n^\epsilon x) \leq d(\bar{h}x, r) + d(r, hg_n^\epsilon x) <$$
$$< (|g_n|_x/2 + d_2 + 2) + (|g_n|_x/2 - (3 + d_2)) = |g_n|_x - 1$$

which contradicts the minimality.                                    $\square$

Recall that by Convention 18 $N$ is an integer $N \geq N_1$.

We now proceed by assigning to a product $w = hg_n^\epsilon$ a segment $S_w$ of length $N$ that lies in the 2-neighborhood of any geodesic segment $[x, wx]$. Note that $S_w$ is defined for the product $w = hg_n^\epsilon$ with $h \in H$ and $\epsilon \in \{-1, 1\}$ and not for the element $w$. These stable parts play an important role in our proof of Proposition 13. We will see that stable parts in some sense "survive" in long products involving $g_n$ and elements of $H$. This will later on allow us to use the "pasting of quasigeodesics" argument to show that $U = \langle M \rangle$ is a free product $U = H * \langle g_n \rangle$ and therefore $U$ is in fact free. In the following definition case 1 is a special case of case 2. However we believe that the overlap creates some additional transparency.

**Definition 20.** [Stable part] Let $N \geq N_1$ be an integer. Suppose that $w$ is the product $w = hg_n^\epsilon$ with $h \in H$ and $\epsilon \in \{-1, 1\}$ and that $|g_n|_x \geq 4N + 2d_2(n-1) + 2d_4(n-1, d_2(n-1) + 11) + 10$. We assign to $w$ the *stable part $S_w$ of $w$ relative $N$* as follows:

1. If $h = 1$ then $w = g_n^\epsilon$. We choose a geodesic segment $[x, g_n^\epsilon x]$ and put
$$S_w := [s, t] \subset [x, g_n^\epsilon x]$$

where $s, t \in [x, g_n^\epsilon x]$ are chosen such that $d(x, s) = |g_n|_x/2 - N - d_2(n - 1)$ and $d(x, t) = |g_n|_x/2 - d_2(n - 1)$. (In Figure 4 $m$ is the midpoint of $[x, g_n^\epsilon x = wx]$ and $d(m, t) = d_2(n - 1)$ and $d(s, t) = N$.)

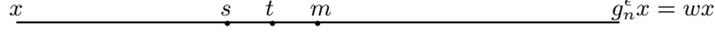

FIGURE 4. Stable part: first case

2. If $h \neq 1$ and the terminal segment of length $|g_n|_x/2 + N + d_2(n - 1)$ of some geodesic segment $h[x, g_n^\epsilon x]$ lies in the 2-neighborhood of some geodesic segment $[x, hg_n^\epsilon x]$, then we choose such a segment $[x, g_n^\epsilon x]$ and put

$$S_w := [s, t] \subset h[x, g_n^\epsilon x]$$

where $[s, t] \subset [x, g_n^\epsilon x]$ are chosen such that $d(hg_n^\epsilon x, s) = |g_n|_x/2 + N + d_2(n - 1)$ and $d(hg_n^\epsilon x, t) = |g_n|_x/2 + d_2(n - 1)$. (In Figure 5 $m$ is the midpoint of $h[x, g_n^\epsilon x]$, $s$ and $t$ are as above.)

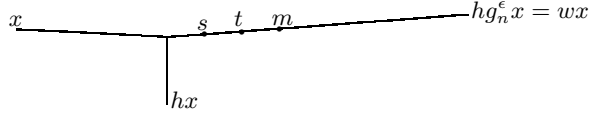

FIGURE 5. Stable part: second case

3. If $h \neq 1$ and the terminal segment of length $|g_n|_x/2 + N + d_2(n - 1)$ of no geodesic segment $h[x, g_n^\epsilon x]$ lies in the 2-neighborhood of any geodesic segment $[x, hg_n^\epsilon x]$ then we choose a path $[x, r] \cup [r, wx]$ corresponding to the product $w = hg_n^\epsilon$ (as in Lemma 15) and put

$$S_w := [s, t] \subset [x, r]$$

where $[s, t] \subset [x, r]$ are chosen such that $d(t, r) = d_4(n - 1, d_2(n - 1) + 11)$ and $d(s, r) = d_4(n - 1, d_2(n - 1) + 11) + N$. (This situation is illustrated in Figure 6.) To ensure the existence of such $[s, t]$ we have to verify that $d(x, r) \geq d_4(n - 1, d_2(n - 1) + 11) + N$. Suppose that $d(x, r) < d_4(n - 1, d_2(n - 1) + 11) + N$. It follows that $d(hx, r) < d_4(n - 1, d_2(n - 1) + 11) + N + 5$ since otherwise $|hg_n^\epsilon|_x < |g_n^\epsilon|_x - 1$. Since further $d(r, hg_n^\epsilon x) \leq |g_n|_x/2 + N + d_2(n - 1)$ by assumption this implies that $|g_n|_x = d(x, g_n^\epsilon x) = d(hx, hg_n^\epsilon x) \leq d(hx, r) + d(r, hg_n^\epsilon x) \leq |g_n|_x/2 + d_2(n - 1) + d_4(n - 1, d_2(n - 1) + 11) + 2N + 5$. Therefore $|g_n|_x < 4N + 2d_2(n - 1) + 2d_4(n - 1, d_2(n - 1) + 11) + 10$ which contradicts our assumption.

We define *stable parts relative $N$*, denoted $S_v$, for a product $v = g_n^\epsilon h$ relative $N$ by considering $v$ as the inverse of the product $w = v^{-1} = g_n^{-\epsilon} h^{-1}$: Since $S_w = [s, t]$ lies in the 2-neighborhood of $[x, wx]$ for some geodesic $[x, wx]$, it follows that $w^{-1}[s, t]$ lies in 2-neighborhood of $[w^{-1}x, x] = [vx, x]$ for some geodesic segment $[vx, x]$. We put $s' := w^{-1}t$, $t' := w^{-1}s$ and define

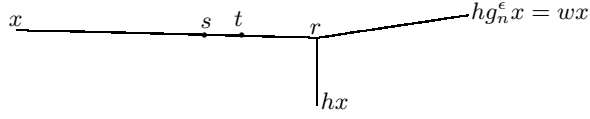

$S_v := w^{-1}[t, s] = [s', t']$, where $[s, t]$ is the geodesic segment $[s, t]$ traveled from $t$ to $s$.

**Remark 21.** Note that a stable part relative $N$ is defined as a geodesic segment of length $N$. It must be stressed that stable parts are defined for ordered products rather than group elements. Thus for $w = 1 \cdot g_n$ and $v = g_n \cdot 1$ the stable parts $S_v$ and $S_w$ are not the same. Indeed, the stable part $S_w$ is a segment of length $N$ somewhat to the left of the midpoint of $[x, g_n x]$. On the other hand $S_v$ is a segment of length $N$ a little to the right of the midpoint of $[x, g_n x]$. It is also very important to note that stable parts are defined only under the minimality assumption on $M$ made in the beginning of this section.

We will use the following notational conventions.

**Convention 22.** Let $u$ be a product of type
$$u = h_1 g_n^{\epsilon_1} h_2 g_n^{\epsilon_2} \cdots g_n^{\epsilon_{l-1}} h_l g_n^{\epsilon_l} h_{l+1},$$
where $\epsilon_i \in \{-1, 1\}$ for $1 \leq i \leq l$, $h_i \in \langle h_1, \dots, h_{n-1} \rangle$ for $1 \leq i \leq l+1$ and $h_i \neq 1$ if $\epsilon_i = -\epsilon_{i+1}$ for $2 \leq i \leq l$. We define $w_i$ to be the product $h_i g_n^{\epsilon_i}$ for $1 \leq i \leq l$ and rewrite $u$ as
$$(h_1 g_n^{\epsilon_1})(h_2 g_n^{\epsilon_2}) \cdots (h_{l-1} g_n^{\epsilon_{l-1}})(h_l g_n^{\epsilon_l})(h_{l+1}) = w_1 w_2 \cdots w_l h_{l+1}$$
and we define $v_i$ to be the product $g_n^{\epsilon_i} h_{i+1}$ and rewrite $u$ as
$$(h_1)(g_n^{\epsilon_1} h_2) \cdots (g_n^{\epsilon_{l-1}} h_l)(g_n^{\epsilon_l} h_{l+1}) = h_1 v_1 v_2 \cdots v_{l-1} v_l.$$
We further define $z_i := h_i g_n^{\epsilon_i} h_{i+1}$ for $1 \leq i \leq l$ and $y_i = g_n^{\epsilon_i} h_{i+1} g_n^{\epsilon_{i+1}}$ for $1 \leq i \leq l-1$.

**Convention 23.** For an integer $N \geq N_1$ put
$$k := 6N + 2d_2(n-1) + d_3(n-1) + 4d_4(n-1, d_2(n-1) + 11) + 35$$
and fix this notation till the end of this section. Note that $k = k(N, n)$ depends on $N$ and $n$. Observe also that $k(N, n)$ is an increasing function of $N$, that is $N \leq N'$ implies $k(N, n) \leq k(N', n)$.

The proof of Proposition 13 relies on the following two lemmas. These lemmas deal with products of length three of the type $h g_n^{\epsilon} h'$ and $g_n^{\epsilon} h g_n^{\delta}$. They basically state that under the minimality assumption on $M$ made in this section the stable parts of $h g_n^{\epsilon}$ and $g_n^{\epsilon} h'$ "survive" and are "disjoint" in $[x, h g_n^{\epsilon} h' x]$ (and a similar statement for $g_n^{\epsilon} h g_n^{\delta}$). Later these facts will allow us to use a "pasting of local quasigeodesics" argument, that is Lemma 6, in the proof of Proposition 13.

**Lemma 24.** *Let $N \geq N_1$ be an integer.*

*Then either $M$ is Nielsen-equivalent to $(g'_1, \dots, g'_n)$ with $||g'_1|| \leq k(N, n)$ (where $k$ is the constant from Convention 23) or for any product $u$ as in Convention 22 the stable parts (relative $N$) are defined for all $v_i$ and $w_i$ and the following holds:*

*There exists a geodesic path $[x, a] \cup [a, z_i x]$ such that $S_{w_i}$ lies in the 15-neighborhood of $[x, a]$ and $h_i S_{v_i}$ lies in the 15-neighborhood of $[a, z_i x]$.*

*Proof.* Recall that in this lemma, as it was assumed throughout this section, $M$ is a minimal tuple with $|g_i|_x \leq |g_n|_x$ for $i = 1, \dots, n-1$ and that Theorem 11 holds for $n-1$.

If $|g_n|_x \leq k$, Lemma 24 obviously holds. Suppose now that $|g_n|_x > k$. Hence the stable parts exists for all $v_i$ and $w_i$ by Definition 20.

We now consider the products $z_i = h_i g_n^{\epsilon_i} h_{i+1}$ and establish the conclusion of Lemma 24. We will first take care of the special case that $h_i$ is trivial (the case that $h_{i+1}$ is trivial is symmetrical). Then we will look at the case that $h_i$ and $h_{i+1}$ are non-trivial and study the two situations which Lemma 16 gives for products of length 3.

Suppose that $h_i$ is trivial. Since in this case $z_i = g_n^{\epsilon_i} h_{i+1} = v_i$, it is clear that $h_i S_{v_i} = S_{v_i}$ lies in the 2-neighborhood of $[x, z_i x]$. The stable part $S_{w_i}$ of $w_i = h_i g_n^{\epsilon_i} = g_n^{\epsilon_i}$ lies by definition in the left part of some geodesic segment $[x, g_n^{\epsilon_i} x]$ in distance $d_2(n-1)$ from the middle. Lemma 19 implies that $S_{w_i}$ lies in the 2-neighborhood of $[x, z_i x]$. Since we also know that $S_{v_i}$ never lies in the left half of $[x, g_n^{\epsilon_i} x]$, the assertion of Lemma 24 holds.

Suppose now that $h_i$ and $h_{i+1}$ are both non-trivial. We apply Lemma 16 for $g = h_i$, $h = g_n^{\epsilon_i}$ and $f = h_{i+1}$ and we choose $p$ and $q$ as in its conclusion.

**Case 1:** Situation 1 of Lemma 16 occurs for the product $h_i g_n^{\epsilon_i} h_{i+1}$, so that we have the picture as in Figure 7.

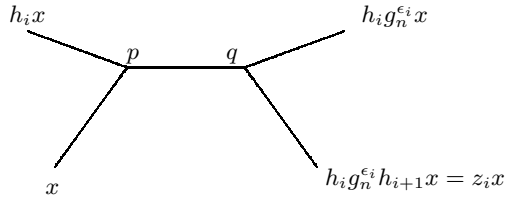

FIGURE 7. Cancellation in $z_i$: first case

We will show that $S_{w_i}$ lies in the 9-neighborhood of any geodesic segment $[x, z_i x]$ (as the result then follows for $h_i S_{v_i}$ by symmetry). By the definition of stable parts and by Lemma 16 we know that $S_{w_i}$ either lies in the 5-neighborhood of $[x, p]$ or in the 5-neighborhood of $[p, q] \cup [q, h_i g_n^{\epsilon_i} x]$. In the first case $S_{w_i}$ lies in the 9-neighborhood of $[x, z_i x]$ since $[x, p]$ lies in the 4-neighborhood of $[x, z_i x]$. Suppose now that $S_{w_i}$ lies in the 5-neighborhood of $[p, q] \cup [q, h_i g_n^{\epsilon_i} x]$. Then $S_{w_i}$ in fact lies in the 9-neighborhood of $[p, q]$, since by Lemma 19 $d(q, h_i g_n^{\epsilon_i} x) \leq |g_n|_x/2 + d_2(n - 1)$ and since by definition the stable part lies on the left half of $h_i[x, g_n^{\epsilon_i} x]$ with distance at least $d_2(n - 1)$ from the midpoint.

To prove the assertion of Lemma 24 we now only have to look at the case when the stable parts $S_{w_i}$ and $h_i S_{v_i}$ both lie in the 2-neighborhood of $h_i[x, g_n^{\epsilon_i} x]$ and therefore in the 9-neighborhood of $[p, q]$. In this situation, however, the desired statement follows as before since the stable parts lie in "different halves" of $h_i[x, g_n^{\epsilon_i} x]$.

**Case 2**: Situation 2 of Lemma 16 occurs for the product $h_i g_n^{\epsilon_i} h_{i+1}$, as shown in Figure 8. This is the only situation where we need the fact that statement (4) of Theorem 11 holds for $n - 1$. Thus this case reveals the reasons why we need the complex statement of Theorem 11 in order to carry out the induction.

It follows from Lemma 19 that $d(q, h_i g_n^{\epsilon_i}) \leq |g_n|_x/2 + d_2(n - 1) + 2$. Since we are in the second case of Lemma 16, this implies that $d(h_i x, q) \geq |g_n|_x/2 - d_2(n - 1) - 2$ and hence the stable part $S_{w_i}$ of $w_i$ is a subsegment of some geodesic segment $[x, h_i x]$. It clearly suffices to show that $S_{w_i}$ lies in the 10-neighborhood of $[x, p]$ (and hence by symmetry $h_i S_{v_i}$ lies in the 10-neighborhood of $[p, z_i x]$). Because of the definition of the stable part it is enough to show that $d(p, q) \leq d_4(n - 1, d_2(n - 1) + 11)$.

Suppose, on the contrary, that $d(p, q) > d_4(n - 1, d_2(n - 1) + 11)$. By Lemma 16 there exists a subsegment $SS$ of $h_i g_n^{\epsilon_i}[x, h_{i+1} x]$ of length $d_4(n - 1, d_2(n - 1) + 11)$ in the 5-neighborhood of $[p, q]$. Since assertion 4 of Theorem 11 holds for $n - 1$, there exists a point $p' \in SS$ and an element $h_i g_n^{\epsilon_i} h' \in h_i g_n^{\epsilon_i} H$ such that

$$d(p', h_i g_n^{\epsilon_i} h' x) \leq \max_{i=1,\ldots,n-1}(|g_i|_x/2 - (d_2(n - 1) + 11) \leq |g_n|_x/2 - d_2(n - 1) - 11.$$

We now choose $\bar{p} \in [p, q]$ such that $d(\bar{p}, p') \leq 5$. It follows that

$$d(\bar{p}, h_i g_n^{\epsilon_i} h' x) \leq |g_n|_x/2 - (d_2(n - 1) + 11) + 5 = |g_n|_x/2 - d_2(n - 1) - 6$$

Lemma 16 also guarantees that there exists a point $p'' \in [x, h_i x]$ such that $d(\bar{p}, p'') \leq 5$. Since assertion 1 of Theorem 11 holds for $n - 1$, there exists an element $\bar{h} \in H$ such that

$$d(p'', \bar{h} x) \leq \max_{i=1,\ldots,n-1}(|g_i|_x/2 + d_2(n - 1)) \leq |g_n|_x/2 + d_2(n - 1).$$

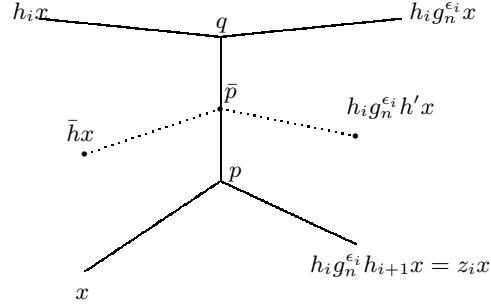

FIGURE 8. Cancellation in $z_i$: second case

This implies that $d(\bar{p}, \bar{h}x) \leq |g_n|_x/2 + d_2(n-1) + 5$. Hence

$$|\bar{h}^{-1}h_ig_n^{\epsilon_i}h'|_x = d(x, \bar{h}^{-1}h_ig_n^{\epsilon_i}h'x) = d(\bar{h}x, h_ig_n^{\epsilon_i}h'x) \leq$$

$$d(\bar{h}x, \bar{p}) + d(\bar{p}, h_ig_n^{\epsilon_i}h'x) \leq$$

$$\leq |g_n|_x/2 + d_2(n-1) + 5 + |g_n|_x/2 - d_2(n-1) - 6 = |g_n|_x - 1,$$

which contradicts the minimality of $M$. $\qquad\square$

**Lemma 25.** *Let $N \geq N_1$ be an integer.*

*Then either $M$ is Nielsen-equivalent to $(g'_1, \ldots, g'_n)$ with $||g'_1|| \leq k(N, n)$ (where $k$ is the constant from Convention 23) or for any product $u$ as in Convention 22 the stable parts (relative $N$) are defined for all $v_i$ and $w_i$ and the following holds.*

*There exists a geodesic path $[x, a] \cup [a, b] \cup [b, y_ix]$ such that:*

*(1) the stable part $S_{v_i}$ lies in the 15-neighborhood of $[x, a]$;*

*(2) a subsegment of length at least*

$$|h_{i+1}|_x - |g_n|_x - 2N - 2d_4(n-1, d_2(n-1) + 11) - 2d_2(n-1)$$

*of $g_n^{\epsilon_i}[x, h_{i+1}x]$ lies in the 15-neighborhood of $[a, b]$;*

*(3) the segment $g_n^{\epsilon_i}S_{w_{i+1}}$ lies in the 15-neighborhood of $[b, y_ix]$.*

*Proof.* Recall that by Convention 14 $M$ is a minimal tuple with $|g_i|_x \leq |g_n|_x$ for $i = 1, \ldots, n-1$ and that Theorem 11 holds for $n-1$.

If $|g_n|_x \leq k$, Lemma 25 obviously holds. Suppose now that $|g_n|_x > k$. Hence the stable parts exists for all $v_i$ and $w_i$ by Definition 20.

We will now look at the products $y_i = g_n^{\epsilon_i}h_{i+1}g_n^{\epsilon_{i+1}}$ and establish the conclusion of Lemma 25. Again we apply Lemma 16 after first taking care of the special case when $h_{i+1}$ is trivial.

If $h_{i+1}$ is trivial, the assertion of Lemma 25 is clear since otherwise $g_n$ has translation length at most $2N + 2d_2(n-1) + 5 \leq k$ by part (1) of Lemma 19.

Suppose that $h_{i+1}$ is non-trivial. We apply Lemma 16 for $g = g_n^{\epsilon_i}$, $h = h_{i+1}$ and $f = g_n^{\epsilon_{i+1}}$ and we choose $p$ and $q$ as in the conclusion of Lemma 16.

**Case 1:** Suppose that we are in the first situation of Lemma 16. If $d(p,q) \geq 2N + 2d_4(n-1, d_2(n-1) + 11)$ then the assertion of Lemma 25 follows immediately from the definition of the stable part. Indeed, in this case the stable parts $S_{v_i}$ and $g_n^{\epsilon_i} S_{w_{i+1}}$ have length $N$ and lie (if they lie in the 4-neighborhood of $g_n^{\epsilon_i}[x, h_{i+1}x]$) in distance $d_4(n-1, d_2(n-1) + 11)$ from $p$ and $q$, respectively.

Suppose now that $d(p,q) < 2N + 2d_4(n-1, d_2(n-1) + 11)$. There are two cases to consider.

**Case 1.A.** Suppose first that $d(p,q) < 2N + 2d_4(n-1, d_2(n-1) + 11)$ and $\epsilon_i = -\epsilon_{i+1}$. (This situation is illustrated in Figure 9.) We will show that in this case $\|h_{i+1}\|$ is small.

Suppose that $d(g_n^{\epsilon_i}x, p) \geq d(g_n^{\epsilon_i}h_{i+1}x, q)$ (as the opposite case is symmetrical). Note that $d(g_n^{\epsilon_i}x, p) \leq |h_{i+1}|_x/2 + 5$, since otherwise $|g_n^{\epsilon_i}h_{i+1}|_x < |g_n^{\epsilon_i}|_x - 1$ which contradicts the minimality of $M$. This implies that $d(p, g_n^{\epsilon_i}h_{i+1}x) \geq |h_{i+1}|_x/2 - 5$ since $d(g_n^{\epsilon_i}x, g_n^{\epsilon_i}h_{i+1}x) = d(x, h_{i+1}x) = |h_{i+1}|_x$.

It follows that

$$d(g_n^{\epsilon_i}h_{i+1}x, q) \geq d(p, g_n^{\epsilon_i}h_{i+1}x) - d(p,q) \geq$$
$$\geq |h_{i+1}|_x/2 - 5 - (2N + 2d_4(n-1, d_2(n-1) + 11)).$$

Define $z := g_n^{\epsilon_i}h_{i+1}^{-1}g_n^{-\epsilon_i}q$, as shown in Figure 9. Note that

$$g_n^{\epsilon_i}h_{i+1}^{-1}g_n^{-\epsilon_i}[g_n^{\epsilon_i}h_{i+1}x, g_n^{\epsilon_i}h_{i+1}g_n^{-\epsilon_i}x] = [g_n^{\epsilon_i}x, x].$$

This implies, in particular, that $z$ lies in the 5-neighborhood of $[x, g_n^{\epsilon_i}x]$ in distance $d(g_n^{\epsilon_i}h_{i+1}x, q)$ from $g_n^{\epsilon_i}x$. Since $d(g_n^{\epsilon_i}x, p) \geq d(g_n^{\epsilon_i}h_{i+1}x, q)$, we know that $z$ actually lies in the 10-neighborhood of $[p, g_n^{\epsilon_i}x]$. Hence there exists a point $\bar{p} \in [p, g_n^{\epsilon_i}x]$ such that $d(z, \bar{p}) \leq 10$ and therefore $d(\bar{p}, g_n^{\epsilon_i}x) \geq d(g_n^{\epsilon_i}x, z) - 10 = d(g_n^{\epsilon_i}h_{i+1}x, q) - 10$. It is clear that

$$d(\bar{p}, q) \leq |h_{i+1}|_x - d(\bar{p}, g_n^{\epsilon_i}x) - d(g_n^{\epsilon_i}h_{i+1}x, q) + 5 \leq |h_{i+1}|_x - 2d(g_n^{\epsilon_i}h_{i+1}x, q) + 15 \leq$$
$$\leq |h_{i+1}|_x - 2(|h_{i+1}|_x/2 - (2N + 2d_4(n-1, d_2(n-1) + 11)) - 5) + 15 =$$
$$= 4N + 4d_4(n-1, d_2(n-1) + 11) + 25,$$

and therefore

$$d(q, g_n^{\epsilon_i}h_{i+1}^{-1}g_n^{-\epsilon_i}q) = d(q,z) \leq 4N + 4d_4(n-1, d_2(n-1) + 11) + 35.$$

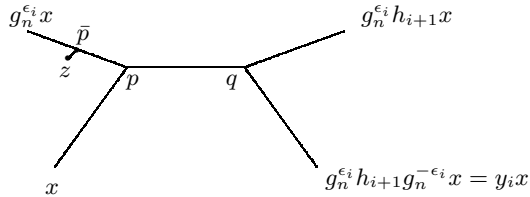

FIGURE 9. Cancellation in $y_i$: case 1.A

This implies that

$$||h_{i+1}|| \leq d(g_n^{-\epsilon_i}q, h_{i+1}g_n^{-\epsilon_i}q) =$$
$$= d(q, g_n^{\epsilon_i}h_{i+1}g_n^{-\epsilon_i}q) \leq 4N + 4d_4(n-1, d_2(n-1) + 11) + 35.$$

Recall that by assumption made in the beginning of this section assertion 3 of Theorem 11 holds for $n-1$ with respect to any base-point of $X$, in particular with respect to the base-point $x' = g_n^{-\epsilon_i}q$. The above inequality implies that $d(x', h_{i+1}x') = |h_{i+1}|_{x'} \leq 4N + 4d_4(n-1, d_2(n-1) + 11) + 35$. Suppose that $M_{n-1} = (g_1, \ldots, g_{n-1})$ is not Nielsen-equivalent to a tuple with the first element of translation length at most $4N + 4d_4(n-1, d_2(n-1) + 11) + 35 + d_3(n-1)$. Choose $j, 1 \leq j \leq n-1$ such that $g_j$ occurs in the freely reduced expression of $h_{i+1}$ as a product of elements of $M_{n-1}$. Then by part 3 of Theorem 11 (applied when $x'$ is chosen as the base-point of $X$) we have $|h_{i+1}|_{x'} \geq |g_j|_{x'} - d_3(n-1)$. Hence

$$|g_j|_{x'} \leq |h_{i+1}|_{x'} + d_3(n-1) \leq 4N + 4d_4(n-1, d_2(n-1) + 11) + 35 + d_3(n-1),$$

which contradicts our assumption on $M_{n-1}$.

Thus $M_{n-1}$ is in fact Nielsen-equivalent to a tuple with the first element of translation length at most $4N + 4d_4(n-1, d_2(n-1)+11) + 35 + d_3(n-1) \leq k$ and hence the conclusion of Lemma 25 holds for $M$ in this case.

**Case 1.B.** Suppose now that $d(p, q) < 2N + 2d_4(n-1, d_2(n-1) + 11)$ and $\epsilon_i = \epsilon_{i+1}$. This case is illustrated in Figure 10.

If $S_{v_i}$ lies in the 5 neighborhood of $[x, p]$ and $g_h^{\epsilon_i}S_{w_{i+1}}$ lies in the 5-neighborhood of $[q, y_i x]$, there is nothing to prove. By definition of stable parts we can assume that either $d(x, p) < |g_n|_x/2 + N + d_2(n-1) + 5$ or that $d(q, y_i x) < |g_n|_x/2 + N + d_2(n-1) + 5$. Without loss of generality we can also assume that $d(q, y_i x) \leq d(x, p)$ and hence $d(q, y_i x) < |g_n|_x/2 + N + d_2(n-1) + 5$. Lemma 19 implies that $d(q, y_i x) \geq |g_n|_x/2 - d_2(n-1) - 5$. Therefore $d(q, m_2) \leq N + d_2(n-1) + 15$ where $m_2$ is the midpoint of a geodesic segment $g_n^{\epsilon_i}h_{i+1}[x, g_n^{\epsilon_{i+1}}x = g_n^{\epsilon_i}x]$.

Hence $d(m_1, z) \leq N + d_2(n-1) + 15$, where $z$ and $m_1$ are defined as $z := h_{i+1}^{-1}g_n^{-\epsilon_{i+1}}q$ and $m_1 := h_{i+1}^{-1}g_n^{-\epsilon_i}m_2$. Note that $m_1$ is approximately the midpoint of $[x, g_n^{\epsilon_i}x]$ (see Figure 10).

Now $d(q, y_i x) \leq d(x, p)$ implies that $d(q, g_n^{\epsilon_i}h_{i+1}x) \geq d(p, g^{\epsilon_i}x) - 5$, that is $d(p, g^{\epsilon_i}x) - d(q, g_n^{\epsilon_i}h_{i+1}x) \leq 5$. Lemma 19 implies that $d(q, g_n^{\epsilon_i}h_{i+1}x) \leq |h_{i+1}|_x/2 + 2$. Since clearly

$$d(p, g^{\epsilon_i}x) + d(p, q) + d(q, g_n^{\epsilon_i}h_{i+1}x) \geq d(g^{\epsilon_i}x, g_n^{\epsilon_i}h_{i+1}x) = |h_{i+1}|_x,$$

we get $d(p, g^{\epsilon_i}x) \geq |h_{i+1}|_x - d(p, q) - d(q, g_n^{\epsilon_i}h_{i+1}x)$.

Thus

$$d(q, g_n^{\epsilon_i}h_{i+1}x) - d(p, g^{\epsilon_i}x) \leq d(p, g^{\epsilon_i}x) - (|h_{i+1}|_x - d(p, q) - d(q, g_n^{\epsilon_i}h_{i+1}x) =$$
$$= 2d(p, g^{\epsilon_i}x) - |h_{i+1}|_x + d(p, q) \leq |h_{i+1}|_x + 4 - |h_{i+1}|_x + d(p, q) = d(p, q) + 4.$$

It follows that $|d(q, g_n^{\epsilon_i} h_{i+1} x) - d(p, g^{\epsilon_i} x)| \leq d(p, q) + 5$. Since $d(m_1, g_n^{\epsilon_i} x) = d(m_2, g_n^{\epsilon_i} h_{i+1} x)$ it follows by symmetry that

$$d(m_1, p) \leq d(q, m_2) + |d(q, g_n^{\epsilon_i} h_{i+1} x) - d(p, g^{\epsilon_i} x)| + 5 \leq$$
$$N + d_2(n-1) + 15 + d(p, q) + 5 + 5 \leq$$
$$\leq N + d_2(n-1) + 2N + 2d_4(n-1, d_2(n-1) + 11) + 25 =$$
$$= 3N + d_2(n-1) + 2d_4(n-1, d_2(n-1) + 11) + 25$$

It is further clear that

$$d(q, z) \leq d(q, p) + d(p, m_1) + d(m_1, z) \leq$$
$$\leq (2N + 2d_4(n-1, d_2(n-1) + 11)) +$$
$$+ (3N + d_2(n-1) + 2d_4(n-1, d_2(n-1) + 11) + 25) + (N + d_2(n-1) + 5) =$$
$$= 6N + 4d_4(n-1, d_2(n-1) + 11) + 2d_2(n-1) + 30 \leq k$$

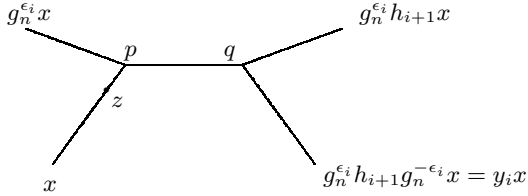

FIGURE 10. Cancellation in $y_i$: case 1.B

Since $z := h_{i+1}^{-1} g_n^{-\epsilon_{i+1}} q$, it follows that $||g_n^{\epsilon_i} h_{i+1}|| = ||h_{i+1}^{-1} g_n^{-\epsilon_{i+1}}|| \leq k$. Thus $M$ is Nielsen-equivalent to a tuple $(g_1', \ldots, g_n')$ with $||g_1'|| \leq k$ and the conclusion of Lemma 25 holds.

**Case 2:** Suppose that we are in the second situation of Lemma 16 for the product $g_n^{\epsilon_i} h_{i+1} g_n^{\epsilon_{i+1}}$, that is that most or all of $h_{i+1}$ "cancels" in $y_i$. We can then argue exactly as in Case 1.A. It turns out that the constants are even smaller since the additive part coming from $d(p, q)$ in the previous case does not occur here. Note that $q$ takes the place of $p$ and $q$ in the argument above. We leave the details of this case to the reader. □

## 5. Proof of Proposition 13

In this section, unless specified otherwise, we assume that $N \geq N_1 > 0$ is a positive integer and that $M$ is a minimal $n$-tuple $M = (g_1, \ldots, g_n)$ not Nielsen equivalent to a tuple containing an element of length at most $k(n, N)$ (the constant provided by Convention 23). Note that by Lemma 24 and Lemma 25 this implies that for any product $u$ as in Convention 22 the stable parts $S_{w_i}$, $S_{v_i}$ relative $N$ are defined for $1 \leq i \leq l$. We will also assume that $|g_1|_x \leq |g_2|_x \leq \ldots |g_n|_x$. As before, we will denote $M_{n-1} = (g_1, \ldots, g_{n-1})$ and $H = \langle M_{n-1} \rangle$.

**Convention 26.** Let $N \geq N_1$ and

$$u = h_1 g_n^{\epsilon_1} h_2 g_n^{\epsilon_2} \cdots g_n^{\epsilon_{l-1}} h_l g_n^{\epsilon_l} h_{l+1}$$

be a product as in Convention 22. Thus $h_i \in H$ for $i = 1, \ldots, l+1$, $\epsilon_i \in \{1, -1\}$ for $i = 1, \ldots, l$ and $h_i \neq 1$ whenever $\epsilon_{i+1} = -\epsilon_i$ for $i = 2, \ldots, l$.

Also let $v_i = g_n^{\epsilon_i} h_{i+1}$, $w_i = h_i g_n^{\epsilon_i}$, $z_i = h_i g_n^{\epsilon_i} h_{i+1}$ and $y_i = g_n^{\epsilon_i} h_{i+1} g_n^{\epsilon_{i+1}}$ be as in Convention 22. For each $i$ let $S_{w_i} = [s_i, t_i]$ and let $S_{v_i} = [s'_i, t'_i]$ be the stable parts relative $N$. Denote $w_1 w_2 \ldots w_{i-1} s_i = S_i$, $w_1 w_2 \ldots w_{i-1} t_i = T_i$, $w_1 \ldots w_{i-1} h_i s'_i = S'_i$ and $w_1 \ldots w_{i-1} h_i t'_i = T'_i$. We define

$$\sigma_N = [x, S_1] \cup [S_1, T_1] \cup [T_1, S'_1] \cup [S'_1, T'_1] \cup [T'_1, S_2] \cup [S_2, T_2] \cup$$
$$\cup [T_2, S'_2] \cup [S'_2, T'_2] \cup \cdots \cup [S'_l, T'_l] \cup [T'_l, ux]$$

The following lemma, despite its technical appearance, is more or less a re-statement of Lemma 24, Lemma 25 in the form convenient for proving Proposition 13.

**Lemma 27.** *Let $N \geq N_1$ and let $\sigma_N$ be as above. Then the following hold:*

1. *The path $\sigma_N$ is a $N$-local $(1, 100)$-quasigeodesic. Therefore, by the choice of $N_1$ in Convention 7, the path $\sigma_N$ is a $(K, K)$-quasigeodesic and it is $L$-Hausdorff close to any geodesic $[x, ux]$.*
2. $d(x, S_1) \geq |g_n|_x / 2 - 2N - 2d_2(n-1) - d_4(n-1, d_2(n-1) + 11) - 24$.
3. $d(x, S_1) \geq |h_1|_x - |g_n|_x / 2 - N - d_2(n-1) - d_4(n-1, d_2(n-1, 11)) - 10$.
4. $d(T'_l, ux) \geq |g_n|_x / 2 - 2N - 2d_2(n-1) - d_4(n-1, d_2(n-1) + 11) - 24$.
5. $d(T'_l, ux) \geq |h_{l+1}|_x - |g_n|_x / 2 - N - d_2(n-1) - d_4(n-1, d_2(n-1, 11)) - 10$.

*Proof.* Parts (1) follows from Lemma 24 and Lemma 25. Statements (2)-(5) follow from the definition of the stable part.  $\square$

We now have all tools necessary to finish off:

*Proof of Proposition 13.* Recall that $L$ and $N_1$ are the constants specified in Convention 7 and that $T > N_1$ by hypothesis.

Suppose that $M = (g_1, \ldots, g_n) \in G^n$ be a minimal $n$-tuple where $|g_i|_x \leq |g_n|_x$ for $i = 1, \ldots, n-1$. Denote $U = \langle M \rangle \leq G$, $M_{n-1} = (g_1, \ldots, g_{n-1})$ and $H = \langle M_{n-1} \rangle \leq U \leq G$. Put $c_1(n, T) := \max\{k(3T, n), d_1(n-1, T)\}$, where $k(3T, n)$ is the constant given by Convention 23. Also put $c_2(n) := d_2(n-1) + 100 + L$ and

$$c_3(n) := 4N_1 + d_3(n-1) + 4d_2(n-1) + 2d_4(n-1, d_2(n-1) + 11) + 49 + 2L + K(2L+2).$$

We will show that Proposition 13 holds with these constants.

Assume that $M$ is not equivalent to a tuple with an element of translation length at most $c_1(n, T)$ (otherwise the statement of Proposition 13 obviously holds). Hence $M_{n-1}$ is not equivalent to a tuple with an element of translation length at most $c_1(n, T)$. Since Theorem 11 holds for $n-1$, the group $H$ is free on $M_{n-1}$ and $H$ and quasi-isometrically embedded in $X$ via the orbit map. Moreover, conditions 2-4 of Proposition 13 hold for any $u \in H$ with the constants $c_1(n, T)$, $c_2(n)$ and $c_3(n)$ specified above.

We first show that $U$ is free on $M$. Since we know that $H$ is free on $M_{n-1}$ it suffices to show that $U = H * \langle g_n \rangle$ and that $g_n$ is of infinite order. It is clear that $g_n$ is of infinite order since otherwise $||g_n|| \leq c_1(1, 0) \leq c_1(n, T)$. Consider an arbitrary element $u = h_1 g_n^{\epsilon_1} h_2 g_n^{\epsilon_2} \ldots g_l^{\epsilon_{l-1}} h_l g_n^{\epsilon_l} h_{l+1}$ as in Convention 26. To see that $U = H * \langle g_n \rangle$ we need to prove that $u \neq 1$.

We will employ the notations used in Convention 26 . Note that for each $w_i$ and $v_i$ the stable parts $S_{w_i}$ and $S_{v_i}$ of length $N_1$ are defined since otherwise by Lemma 24 and Lemma 25 the tuple $M$ is Nielsen-equivalent to a tuple with an element of length at most $k(N_1, n) \leq k(3T, n) \leq c_1(n, T)$.

Consider the path $\sigma_{N_1}$ from $x$ to $ux$ defined as in Convention 26.

By Lemma 27 the path $\sigma_{N_1}$ is $N_1$-local $(1, 100)$-quasigeodesic and $(K, K)$-global quasigeodesic and $\sigma_{N_1}$ is $L$-Hausdorff close to $[x, ux]$.

Note that the length of $\sigma_{N_1}$ is at least $N_1$. Therefore $d(x, ux) = |u|_x \geq N_1/K - K > 0$ by the choice of $N_1$ in Convention 7. Thus $u \neq 1$ and so $U$ is free on $M$ as required.

We will now establish part 2 of Proposition 13.

Recall that by assumption condition 2 of Theorem 11 holds for any $h \in H = \langle M_{n-1} \rangle$ with the constant $d_2(n-1)$. Suppose $u \in U, u \notin H$, that is $u$ is an alternating product as in Convention 26, involving at least one $g_n$. We already know that $\sigma_{N_1}$ is $(K, K)$-quasigeodesic and hence $L$-Hausdorff close to $[x, ux]$. Thus $[x, ux]$ lies in $L$-neighborhood of the path $\sigma_{N_1}$. By Lemma 27 $\sigma_{N_1}$ is contained in 100-neighborhood of the path

$$\sigma' = [x, h_1 x] \cup h_1 [x, g_1^{\epsilon_1} x] \cup h_1 g_1^{\epsilon_1} [x, h_2 x] \cup h_1 g_1^{\epsilon_1} h_2 [x, g_2^{\epsilon_2} x] \cup \ldots$$
$$\cdots \cup h_1 g_n^{\epsilon_1} h_2 g_n^{\epsilon_2} \ldots g_l^{\epsilon_{l-1}} h_l [x, g_n^{\epsilon_l} x] \cup h_1 g_n^{\epsilon_1} h_2 g_n^{\epsilon_2} \ldots h_l g_n^{\epsilon_l} [x, h_{l+1} x].$$

Therefore $[x, ux]$ is contained in the $(100 + L)$-neighborhood of $\sigma'$. As we noticed before, condition 2 of Theorem 11 holds for any $h \in H = \langle M_{n-1} \rangle$ with the constant $d_2(n-1)$. Therefore condition 2 of Proposition 13 holds for $u$ with the constant $c_2(n) = d_2(n-1) + 100 + L$ specified as above. Since the orbit $Ux$ is $U$-invariant, this means that the orbit $Ux$ is $a$-quasiconvex in $X$ where $a = |g_n|_x/2 + c_2(n)$. Thus condition 2 of Proposition 13 is verified.

We will now verify that part 3 of Proposition 13 holds with $c_3(n)$ as specified above.

Recall that $|g_n|_x \geq |g_j|_x$ for $j = 1, \ldots, n - 1$. If $u$ is a freely reduced word in $M$ which does not involve $g_n$, then statement 3 of Proposition 13 holds for $u$ with the constant $d_3(n-1) \leq c_3(n)$. Suppose now that $u$ is an alternating product involving $g_n$, as in Convention 26. Again consider the path $\sigma_{N_1}$ from $x$ to $ux$ defined as in Convention 26.

Thus $[S_i, T_i]$ and $[S_i', T_i']$ are stable parts of appropriate elements relative $N_1$ and therefore have length $N_1$. Recall that by Lemma 27 $\sigma_{N_1}$ is a $(K, K)$-quasigeodesic from $x$ to $ux$ which is $L$-close to $[x, ux]$.

By Lemma 27 we know that $d(x, S_1) \geq |g_n|_x/2 - 2N_1 - 2d_2(n-1) - d_4(n-1, d_2(n-1) + 11) - 24$ and $d(T_l', ux) \geq |g_n|_x/2 - 2N_1 - 2d_2(n-1) - d_4(n-$ $1, d_2(n-1) + 11) - 24$. Since $\sigma_{N_1}$ is $L$-close to $[x, ux]$, there are points $p$ and $q$ on $[x, ux]$ such that $d(S_1, p) \leq L$ and $d(T'_l, q) \leq L$. Hence by the triangle inequality we have

$$d(x, p) \geq |g_n|_x / 2 - 2N_1 - 2d_2(n-1) - d_4(n-1, d_2(n-1) + 11) - 24 - L$$
$$d(q, ux) \geq |g_n|_x / 2 - 2N_1 - 2d_2(n-1) - d_4(n-1, d_2(n-1) + 11) - 24 - L$$

Note that by definition of $\sigma_{N_1}$ the paths $[x, S_1]$ and $[T'_l, ux]$ do not overlap in $\sigma$. If $d(x, p) \leq d(x, q)$ then $d(x, ux) \geq d(x, p) + d(q, ux)$ and hence

$$|u|_x = d(x, ux) \geq$$
$$\geq |g_n|_x - 4N_1 - 4d_2(n-1) - 2d_4(n-1, d_2(n-1) + 11) - 48 - 2L \geq$$
$$\geq |g_n|_x - c_3(n),$$

as required.

Suppose that $d(x, p) \geq d(x, q)$, so that $q \in [x, p]$. Since $[x, S_1]$ is a geodesic and $d(p, S_1) \leq L$, there is a point $q'$ on $[x, S_1]$ such that $d(q, q') \leq L + 1$. Hence $d(q', T_l) \leq d(q', q) + d(q, T_l) \leq 2L + 1$. Since $\sigma_{N_1}$ is $(K, K)$-quasigeodesic, the length of the $\sigma_{N_1}$-segment from $q'$ to $T_l$ is at most $K(2L + 1) + K$. The point $S_1$ lies on the $\sigma_{N_1}$-segment from $q'$ to $T'_l$ and therefore $d(q', S_1) \leq K(2L + 1) + K$. Hence

$$d(q, p) \leq d(q, q') + d(q', S_1) + d(S_1, p) \leq$$
$$\leq (L + 1) + [K(2L + 1) + K] + L = K(2L + 2) + 2L + 1$$

Therefore we can estimate the length of $[x, ux]$ as follows:

$$|u|_x = d(x, ux) = d(x, p) + d(q, ux) - d(p, q) \geq$$
$$|g_n|_x - 4N_1 - 4d_2(n-1) - 2d_4(n-1, d_2(n-1) + 11) - 49 - 2L - K(2L + 2) \geq$$
$$\geq |g_n|_x - c_3(n),$$

as required. Thus part 3 of Proposition 13 is verified.

We will now establish part 4 of Proposition 13.

Suppose $u = h_1 g_n^{\epsilon_1} h_2 g_n^{\epsilon_2} \dots g_l^{\epsilon_{l-1}} h_l g_n^{\epsilon_l} h_{l+1}$ is an alternating product as in Convention 26. Also let $\sigma = \sigma_{3T}$ be defined as in Convention 26. Since $T > N_1$, Lemma 27 implies that $\sigma$ is a $(K, K)$-quasigeodesic which is $L$-close to $[x, ux]$.

Suppose $I = [t, r]$ is a subsegment of $[x, ux]$ of length $10T$. Then there are points $t', r'$ on $\sigma$ such that $d(t, t') \leq L$ and $d(r, r') \leq L$. Since $X$ is 1-hyperbolic and geodesic quadrilaterals are 2-thin, the segment $I = [t, r]$ and a geodesic $[t', r']$ are $L + 2$-Hausdorff close. Since the segment of $\sigma$ from $t'$ to $r'$ is a $(K, K)$-quasigeodesic, it is $L$-close to $[t', r']$. Thus the subpath of $\sigma$ from $t'$ to $r'$ is $2L + 2$-Hausdorff close to $I = [t, r]$.

By Convention 7, we have $N_1 \geq 2L$. Since we chose $T > N_1$ in the beginning of the proof of Proposition 13, we have $T \geq 2L$. Hence $d(t', r') \geq 10T - 2L \geq 9T$ and so the length of the part of $\sigma$ from $t'$ to $r'$ is at least $9T$. Hence the part of $\sigma$ from $t'$ to $r'$ either contains a subpath $\tau$ of length $3T$ such that either $\tau = [S_i, T_i]$ or $\tau = [S'_i, T'_i]$ (which are by definition of length $3T$) or $\tau$ is contained in one of $[x, S_1]$, $[T_i, S'_i]$, $[T'_i, S_{i+1}]$ or $[T'_l, ux]$.

It now follows from Lemma 24 and Lemma 25 that the 50-neighborhood of $\tau$ contains a subsegment $\tau'$ of either $w_1 \cdots w_{i-1}[x, z_i x]$ or $w_1 \cdots w_{i-1} h_i[x, y_i x]$ of length $3T$. By Corollary 17 the 24-neighborhood of $\tau'$ contains a subsegment $\alpha$ of length at least $T - 40$ contained in one of one of the segments $h_1 g_n^{\epsilon_1} h_2 g_n^{\epsilon_2} \ldots g_n^{\epsilon_{i-1}} h_i [x, g_n^{\epsilon_i} x]$ or $h_1 g_n^{\epsilon_1} h_2 g_n^{\epsilon_2} \ldots g_n^{\epsilon_{i-1}} [x, h_i x]$. It follows that the 44-neighborhood of $\tau'$ contains such a segment of length at least $T$. Thus we see that the $2L + 96$-neighborhood of $I$ contains $\alpha$ and the assertion (4) of Proposition 13 is established.

Finally, we establish that the orbit map $U \longrightarrow X$, $u \mapsto ux$ is a quasi-isometric embedding. We already know that $U$ is free on $M$. For any $u \in U$ denote by $|u|_U$ the word-length of $u$ with respect to the free basis $M$. Similarly, for $h \in H$ we denote by $|h|_H$ the word-length of $h$ with respect to the free basis $M_{n-1} = (g_1, \ldots, g_{n-1})$ of $H$. Note that since $U$ is free on $M$, we have $|h|_H = |h|_U$ for any $h \in H$.

Since Theorem 11 holds for $n - 1$, we already know that the orbit map for $H$ is a quasi-isometric embedding. Moreover, the $U$-orbit map $x \mapsto ux$ is injective. Indeed, if $u \in U - H$ then $d(x, ux) > 0$ as we have shown when proving that $U$ is free. Suppose $u = h \in H$, $u \neq 1$ fixes $x$, that is $hx = x$. Then $d(h^m x, x) = 0$ for any $m \in \mathbb{Z}$ which contradicts the fact that $H$ is free and the orbit map of $H$ is a quasi-isometric embedding.

Since $H$ with the free group metric is a uniformly discrete space and the orbit map for $H$ is a quasi-isometric embedding, it follows that the $H$-orbit map $h \mapsto hx$ is a bi-Lipschitz bijection onto its image. Thus there is some $C > 1$ be such that for any $h \in H$ we have $|h|_H \leq C|h|_x$.

Since the set $Ux$ is $U$-invariant and the orbit map of $U$ is obviously Lipschitz, it suffices to show that there is a constant $C' > 0$ such that for any $u \in U$ we have $|u|_U \leq C'|u|_x + C'$. If $u \in H$, then obviously $|h|_H = |h|_U \leq C|h|_x$.

Suppose now that $u \notin H$, $u \in U$. Again write down $u$ as a product as in Convention 26. As in the proof that $U$ is free, let $N = N_1$ and let $\sigma = \sigma_N$ be the path from $x$ to $ux$ defined as in Convention 26.

Denote $D := 2N + 2d_4(n - 1, d_2(n - 1) + 11) + 2d_2(n - 1) + 20$. By Lemma 27 if $|h_{i+1}|_x \geq 2[|g_n|_x + D]$ then

$$d(T_i, S'_{i+1}) \geq |h_{i+1}|_x - |g_n|_x - D \geq |h_{i+1}|_x/2.$$

Similarly by Lemma 27, if $|h_1|_x \geq 2[|g_n|_x + D]$ then

$$d(x, S_1) \geq |h_1|_x - |g_n|_x/2 - D \geq |h_1|_x/2.$$

and if $|h_{l+1}|_x \geq 2[|g_n|_x + D]$ then

$$d(T'_l, ux) \geq |h_{l+1}|_x - |g_n|_x/2 - D \geq |h_{l+1}|_x/2.$$

Since the orbit map for $H$ is a quasi-isometry, there are only finitely many elements $h \in H$ such that $|h|_x \leq 2[|g_n|_x + D]$. Let $K'$ be the maximal $H$-length of all such $h \in H$ (so that for every such $h$ we have $|h|_x \leq CK'$). Recall that $\sigma$ is a (K,K)-quasigeodesic. Then for $i = 2, \ldots, l$

either $|h_i|_H \leq K$, $1 \geq |h_i|_H / K'$ or $d(T'_{i-1}, S_i) \geq |h_i|_x / 2 \geq |h_i|_H / (2C)$.

Similarly

either $|h_1|_H \leq K$, $1 \geq |h_1|_H / K'$ or $d(x, S_1) \geq |h_1|_x / 2 \geq |h_1|_H / (2C)$

and

either $|h_{l+1}|_H \leq K$, $1 \geq |h_{l+1}|_H / K'$ or $d(T'_l, ux) \geq |h_{l+1}|_x / 2 \geq |h_{l+1}|_H / (2C)$.

Thus we see that either $|h_i|_H$ is small or a substantial portion of $|h_i|_x$ is "reflected" in a subpath of $\sigma$. This easily implies that the length of $\sigma$ can be estimated from below in terms of $|u|_U$.

Indeed, observe that $|u|_U = l + \sum_{i=1}^{l+1} |h_i|_H$ and that the length $l(\sigma)$ of $\sigma$ can be estimated as follows:

$$l(\sigma) \geq Nl + \sum_{i=2}^{l} d(T'_{i-1}, S_i) + d(x, S_1) + d(T'_l, ux) =$$

$$\geq (N-1)l + \frac{1}{2K'C} \sum_{i=1}^{l+1} |h_i|_H \geq \frac{1}{2K'C} |u|_U$$

Since $\sigma$ is a $(K, K)$-quasigeodesic in $X$, we have

$$|u|_U \leq 2K'Cl(\sigma) \leq 2KK'C|u|_x + K$$

and hence the orbit map for $U$ is indeed a quasi-isometric embedding. This completes the proof of Proposition 13. $\qquad \square$

DEPT. OF MATHEMATICS, UNIVERSITY OF ILLINOIS, 1409 WEST GREEN STREET, URBANA, IL 61801, USA
  *E-mail address*: `kapovich@math.uiuc.edu`

FAKULTÄT FÜR MATHEMATIK, RUHR-UNIVERSITÄT BOCHUM, 44780 BOCHUM, GERMANY
  *E-mail address*: `Richard.Weidmann@ruhr-uni-bochum.de`